\documentclass{svjour3}
\smartqed
\journalname{Probability Theory and Related Fields}
\usepackage{mathptmx}

\usepackage{ifpdf}
\ifpdf
	\usepackage[pdftex]{graphicx}
	\usepackage[pdftex]{hyperref}
	\DeclareGraphicsExtensions{.png,.pdf}
\else
   \usepackage{graphicx}
   \usepackage{hyperref}
   \DeclareGraphicsExtensions{.ps,.eps}
\fi

\usepackage[utf8]{inputenc}
\usepackage{amsmath}
\usepackage{amsfonts}
\usepackage{amssymb}

\usepackage{mathrsfs}
\newcommand{\calA}{\mathscr{A}}
\newcommand{\calB}{\mathscr{B}}

\newcommand{\calE}{\mathscr{E}}

\newcommand{\calG}{\mathscr{G}}

\newcommand{\calR}{\mathscr{R}}

\newcommand{\bbN}{\mathbb{N}}

\newcommand{\bbP}{\mathbb{P}}

\newcommand{\bbR}{\mathbb{R}}

\newcommand{\bbT}{\mathbb{T}}

\newcommand{\bbZ}{\mathbb{Z}}

\newcommand{\bk}[1]{\langle{#1\rangle}{}}
\newcommand{\Zd}{\bbZ^d}

\newcommand{\PF}{\mathbf{Z}}
\newcommand{\IF}[1]{\mathbf{1}_{\{#1\}}}
\newcommand{\setof}[2]{\left\{#1 \,:\, #2 \right\}}
\newcommand{\given}{\,|\,}
\newcommand{\bgiven}{\bigm|}
\newcommand{\dd}{\mathrm{d}}
\newcommand{\etac}{\eta_{\mathrm{c}}}

\newcommand{\normI}[1]{\|#1\|_{\scriptscriptstyle 1}}
\newcommand{\normII}[1]{\|#1\|_{\scriptscriptstyle 2}}
\newcommand{\normINF}[1]{\|#1\|_{\scriptscriptstyle \infty}}
\newcommand{\per}{\mathrm{per}}
\newcommand{\comp}{{\rm c}}
\newcommand{\pinned}{\zeta}

\begin{document}

\title{Wetting of gradient fields: pathwise estimates\thanks{Supported in part by Fonds National Suisse.}}
\titlerunning{Wetting of gradient fields: pathwise estimates}
\author{Yvan Velenik}
\institute{Y. Velenik \at
              Section de Mathématiques, Université de Genève \\
              2-4, rue du Lièvre, Case postale 64\\
              CH-1211 Genève 4\\
              \email{Yvan.Velenik@math.unige.ch}
}
\date{Received: date / Accepted: date}

\maketitle
\begin{abstract}
We consider the wetting transition in the framework of an effective interface model of gradient type, in dimension $2$ and higher. We prove pathwise estimates showing that the interface is localized in the whole thermodynamically-defined partial wetting regime considered in earlier works. Moreover, we study how the interface delocalizes as the wetting transition is approached. Our main tool is reflection positivity in the form of the chessboard estimate.
\end{abstract}

\section{Introduction and results}
\subsection{The model}
Effective interface models of gradient type have been a very active field of research in recent years. In particular, the understanding of the interaction of an interface with various types of external potentials (wall, pinning potential, etc.) has motivated numerous works, resulting in substantial progress on such issues. We refer to~\cite{Fu2005,Gi2001,Ve2006} for reviews of the problems investigated and references.

Among such questions, the analysis of the effect of an attractive wall on the behavior of an interface is of particular relevance. Such a situation is commonly modeled as follows. Let $\Lambda_L = \{-\lceil L/2\rceil+1,\ldots,\lfloor L/2\rfloor\}^d$, $\overline\Lambda_L=\{-\lceil L/2\rceil,\ldots,\lfloor L/2\rfloor+1\}^d$ and $\partial\Lambda_L=\overline\Lambda_L\setminus\Lambda_L$. Interface configurations are given by $\varphi=\{\varphi_i\}_{i\in\Lambda_L}\in \bbR^{\overline\Lambda_L}$. Let also $V:\bbR\to\bbR$ be an even, convex function, with $V\not\equiv0$ and $V(0)=0$. Given $\eta,\lambda\geq 0$, we introduce the following probability measure on $\bbR^{\overline\Lambda_L}$.
\begin{multline*}
\mu^{0}_{L;\lambda,\eta} (\dd\varphi) = \frac1{\PF^{0}_{L;\lambda,\eta}} \exp \bigl[-\frac1{8d}\sum_{\substack{i,j\in\Lambda_L\\i\sim j}} (\varphi_i-\varphi_j)^2 - \lambda \sum_{i\in\Lambda_L} V(\varphi_i) \bigr]\\
\times \prod_{i\in\Lambda_L} \bigl( \dd\varphi_i + \eta\delta_0(\dd\varphi_i) \bigr) \prod_{i\in\partial\Lambda_L} \delta_0(\dd\varphi_i),
\end{multline*}
where $\dd\varphi_i$ and $\delta_0(\dd\varphi_i)$ denote respectively Lebesgue measure and the Dirac mass at $0$, and $i\sim j$ means that $\normI{i-j}=1$. Writing $\Omega_+ = \{\varphi_i\geq 0, \forall i\in\Lambda_L\}$, we then introduce the probability measure
$$
\mu^{+,0}_{L;\lambda,\eta} (\,\cdot\,) = \mu^{0}_{L;\lambda,\eta} (\,\cdot\given \Omega_+).
$$
This is the measure we shall be mostly interested in this work. We shall denote by
\begin{equation}
\label{equ_GibbsMeasure}
\PF^{+,0}_{L;\lambda,\eta} = \mu^{0}_{L;\lambda,\eta}(\Omega_+)
\end{equation}
the associated partition function. Before going on, let us briefly describe the physical meaning of all the pieces entering the definition of $\mu^{+,0}_{L;\lambda,\eta}$. Interpreting as usual $\varphi_i$ as the height of the interface above site $i$, the positivity constraint $\Omega_+$ corresponds to the presence of a hard wall at height $0$, which the interface cannot cross. The term
$$
\frac1{8d}\sum_{\substack{i,j\in\Lambda_L\\i\sim j}} (\varphi_i-\varphi_j)^2
$$
represents the internal energy associated to deformation of the interface from the horizontal plane. The term
$$
\lambda \sum_{i\in\Lambda_L} V(\varphi_i)
$$
represents the contribution to the energy coming from the presence of an external potential. A common choice if $V(x)=x^2$ (usually termed a \textit{mass term}), but given the situation we want to model here a more natural choice is $V(x)=|x|$. The latter choice allows for the interpretation of the interface as separating a thermodynamically stable phase (above) from a thermodynamically unstable phase (below), the latter being stabilized locally because it is favored by the wall; $\lambda$ then measures the difference of free energies between the stable and unstable phases (both being stable when $\lambda=0$); see~\cite{Ve2004,Ve2006} for a more detailed explanation. Finally, for $\eta>0$, the measure
$$
\prod_{i\in\Lambda_L} \bigl( \dd\varphi_i + \eta\delta_0(\dd\varphi_i) \bigr)
$$
models the local attractivity of the interface/wall interaction, by rewarding each contact between the interface and the wall. One way to see this better (which also turns out to be technically useful later) is to realize that $\mu^{0}_{L;\lambda,\eta}$ can be seen as the weak limit of the measures
\begin{multline*}
\mu^{0,(\epsilon)}_{L;\lambda,\eta} (\dd\varphi) = \frac1{\PF^{0,(\epsilon)}_{L;\lambda,\eta}} \exp \bigl[-\frac1{8d}\sum_{\substack{i,j\in\Lambda_L\\i\sim j}} (\varphi_i-\varphi_j)^2 - \lambda \sum_{i\in\Lambda_L} V(\varphi_i) - \sum_{i\in\Lambda_L} U^{(\epsilon)}_\eta(\varphi_i) \bigr]\\
\times \prod_{i\in\Lambda_L} \dd\varphi_i \prod_{i\in\partial\Lambda_L} \delta_0(\dd\varphi_i),
\end{multline*}
where $e^{-U^{(\epsilon)}_\eta(\varphi_i)} = 1+\frac{\eta}{2\epsilon}\IF{|\varphi_i|\leq \epsilon}$, as $\epsilon\downarrow 0$. Similarly, $\mu^{+,0}_{L;\lambda,\eta}$ is easily seen to be given by the weak limit, as $\epsilon\downarrow 0$, of
\begin{equation}
\label{eq_SWpinning}
\mu^{+,0,(\epsilon)}_{L;\lambda,\eta} = \mu^{0,(\epsilon)}_{L;\lambda,2\eta}( \,\cdot\given\Omega_+).
\end{equation}

\subsection{Earlier results}
Various aspects of this model have been studied in several papers. Let us briefly review earlier works relevant to the present contribution. Many of the results quoted below are valid in the more general context of gradient field with uniformly strictly convex interactions, \textit{i.e.} those for which the term $(\varphi_i-\varphi_j)^2$ in the definition of the measure is replaced by $U(\varphi_i-\varphi_j)$ with $U:\bbR\to\bbR$ an even function with second derivative uniformly bounded away from $0$ and $\infty$. To keep the discussion short, we shall not discuss this here (nor shall we discuss the case of non-nearest-neighbor interactions), and refer to the cited papers, and to the reviews mentioned at the beginning, for more information. Let us just remark that most of our analysis actually extends to this case as well, the Gaussian character of the measure being used in an essential way only in very few places. However, most earlier results about the free energy in the wetting problem, upon which our whole approach rests, concern exclusively the Gaussian setting (or Lipschitz interactions).

\subsubsection{Free interface}
We very briefly recall what is known when $\lambda=\eta=0$, for the measure without the positivity constraint, \textit{i.e.}, for the measure $\mu^0_{L;0,0}$. In that case, the measure is Gaussian, and therefore amenable to explicit computations. Many things are known, but for our purposes here, it is enough to say that the variance of the field satisfies\footnote{We write $o_\ell(1)$ to denote a function such that $\lim_{\ell\to\infty}o_\ell(1)=0$.}
$$
\bk{\varphi_0^2}^0_{L;0,0} =
\begin{cases}
(g(1)+o_L(1)) L		&	(d=1),\\
(g(2)+o_L(1)) \log L	&	(d=2),\\
g(d)+o_L(1)		&	(d\geq 3),
\end{cases}
$$
for explicit constants $g(d)>0$, which shows that this measure describes a delocalized interface, with unbounded fluctuations, in dimensions $1$ and $2$, and a localized interface in dimension $3$ and higher. In the latter case, although localized, the interface is strongly correlated,
$$
\lim_{L\to\infty}\bk{\varphi_i\varphi_j}^0_{L;0,0} = (a(d)+o_{\normII{i-j}}(1))\,\normII{i-j}^{2-d},
$$
with $a(d)>0$, $d\geq 3$.
\subsubsection{Interface and pinning potential}
Setting $\eta>0$, keeping everything as before, changes dramatically the behavior of the field \emph{however small $\eta$ is chosen}. More precisely, it is known that the interface is localized in any dimension~\cite{DuMaRiRo1992,BoBr2001,DeVe2000}, and has exponentially decaying covariances~\cite{BoBr2001,IoVe2000}. Moreover, detailed information on the critical behavior as $\eta\downarrow 0$ is available~\cite{BoVe2001}, showing for example that
$$
\lim_{L\to\infty}\bk{\varphi_0^2}^0_{L;0,\eta} =
\begin{cases}
\tfrac12 \eta^{-2}+o(\eta^{-2})			&	(d=1),\\
\frac1\pi |\log\eta| + O(\log|\log\eta|)	&	(d=2),
\end{cases}
$$
and that the rate $m(\eta)$ of exponential decay of $\lim_{L\to\infty}\bk{\varphi_i\varphi_j}^0_{L;0,0}$ satisfies
$$
m(\eta) =
\begin{cases}
\tfrac12 \eta^2	+ o(\eta^2)		&	(d=1),\\
O(\eta^{1/2}/|\log\eta|^{3/4})		&	(d=2),\\
O(\eta^{1/2})				&	(d\geq 3).
\end{cases}
$$
\subsubsection{Interface and hard-wall}
The measure with hard-wall constraint, but no external potentials, \textit{i.e.} $\mu^{+,0}_{L;0,0}$ has been the subject of numerous works, focusing on the associated entropic repulsion phenomenon. Among the results that have been obtained, we highlight the two most relevant in the present context. Let $d\geq 3$; then~\cite{BoDeZe1995,De1996}
$$
\lim_{L\to\infty} \bigl|\frac{\bk{\varphi_0}^{+,0}_{L;0,0}}{\sqrt{\log L}} -2\sqrt{g(d)} \bigr| = 0.
$$
The corresponding result in dimension $2$, whose proof is substantially more intricate, is proved in~\cite{BoDeGi2001} and takes the form
$$
\lim_{L\to\infty} \bigl|\frac{\bk{\varphi_0}^{+,0}_{L;0,0}}{\log L} -2\sqrt{g(2)} \bigr| = 0.
$$
(Actually, the statement in~\cite{BoDeGi2001} has only been proved when the positivity constraint acts on the sub-box $\Lambda_{\delta L}$, $0<\delta<1$, but it is clear that the previous result is true, and that it should be provable in the same way, with some additional, but minor, complications.)

The main thing to observe here is the fact that the interface is repelled by the wall, at a distance that is much larger than its typical fluctuations (which are of order $\sqrt{\log L}$ when $d=2$, and of order $1$ when $d\geq 3$). This is the phenomenon of entropic repulsion. Of course, this does not happen when $d=1$, since the pinned random walk conditioned to be positive converges under diffusive scaling to the Brownian excursion.

\subsubsection{Interface and attractive hard-wall: wetting transition}
We want to describe the behavior of the field when both a hard-wall and a pinning potential are present, $\mu^{+,0}_{L;0,\eta}$. In this situation, there is a competition between the entropic repulsion due to the hard-wall constraint and the localizing effect of the pinning potential.

Let us introduce the finite-volume average density of pinned sites
$$
\rho_L(\eta) = \bk{|\Lambda_L|^{-1} \sum_{i\in\Lambda_L} \IF{\varphi_i=0}}_{L;0,\eta}^{+,0},
$$
and its limit $\rho(\eta)=\lim_{L\to\infty} \rho_L(\eta)$. It is easy to show that $\rho$ is non-decreasing in $\eta$, so the following critical value is well-defined,
$$
\etac = \inf\setof{\eta}{\rho(\eta)>0}.
$$
This critical point can be given an equivalent definition (the equivalence is proved, \textit{e.g.}, in~\cite{CaVe2000}). Let us introduce the free energy (or surface tension, or wall free energy)
$$
f_L(\lambda,\eta) = |\Lambda_L|^{-1}\, \log\frac{\PF^{+,0}_{L;\lambda,\eta}}{\PF^{+,0}_{L;\lambda,0}},
$$
and $f(\lambda,\eta) = \lim_{L\to\infty} f_L(\lambda,\eta)$. Then
$$
\etac = \inf\setof{\eta}{f(0,\eta)>0}.
$$
The sets $\{\eta \leq\etac \}$, resp. $\{\eta>\eta_c\}$, are called regimes of \textit{complete wetting}, resp. \textit{partial wetting}. They are supposed to correspond to regimes in which the interface is delocalized, resp. localized. The phase transition taking place at $\etac$ is known as the \textit{wetting transition}. It is known that $\etac = 0$ when $d\geq 3$~\cite{BoDeZe2000}, while $\etac>0$ when $d=2$~\cite{CaVe2000}\footnote{Actually, it is interesting to observe that it is also proved in~\cite{CaVe2000} that $\etac>0$ in any dimensions if the interaction term $(\varphi_i-\varphi_j)^2$ is replaced by, say, $|\varphi_i-\varphi_j|$.}. The fact that $\etac>0$ in dimension $1$ is easily checked, and has been proved long ago by physicists.

Contrarily to the results described above, there are only very few pathwise results in this setting, except in dimension $1$, where specific features (in particular, a natural renewal structure) makes it possible to fully describe the process~\cite{DeGiZa2005}. Before the present paper, the only pathwise results available are those of~\cite{Ve2004}, and state that, in dimension $2$,
\begin{itemize}
\item For all $\eta$ sufficiently large, the interface is localized, and covariances decay exponentially.
\item For all $\eta <\etac$, the interface delocalizes, in the sense that $\lim_{L\to\infty}\bk{\varphi_0}^{+,0}_{L;0,\eta}=\infty$.
\item For $\eta$ sufficiently small\footnote{$a\asymp b$ meaning here and in the rest of this paper, that there exists a constant $c>0$, depending on nothing except possibly the dimension, such that $ac \leq b \leq a/c$.}, $\bk{\varphi_0}^{+,0}_{L;0,\eta} \asymp\log L$.
\end{itemize}
Notice that $\eta_c=0$ in dimensions $3$ and higher, and thus the analogue of the last statement reduces to the entropic repulsion estimate given above.

The main goal of the present paper is to provide detailed pathwise information on the localized regime in the whole partial wetting regime and in any dimensions. Moreover, we shall give some information on the rate of divergence of the height as the wetting transition is approached (implying among other things that divergence does occur also in dimension $2$).

\subsubsection{Interface, attractive hard-wall, away from coexistence: prewetting}
Finally, letting $\lambda>0$ introduces another source of localization of the interface. Of course, under our assumptions on this potential, it is not surprising that it always yields localization of the interface, which corresponds to the impossibility of growing a large film of thermodynamically unstable phase. The main question here is to understand what happens as the system is brought close to phase coexistence, \textit{i.e.} when $\lambda\downarrow 0$.

The situation studied in~\cite{HrVe2004,Ve2004} is the following: Fix $0\leq\eta<\etac$, in dimension $1$ or $2$, or take $\eta=0$ in dimension $3$ and larger. Set also $\lambda>0$. Then, as $\lambda\downarrow 0$, the system gets closer and closer to the regime of phase coexistence, and in that regime, because of the choice for $\eta$, the interface is delocalized. The problem was then to determine the rate at which this delocalization takes place. The main result of~\cite{Ve2004} can be stated as follows: For all $\lambda>0$ sufficiently small,
$$
\lim_{L\to\infty} \bk{\varphi_0}^{+,0}_{L;\lambda,\eta} \asymp
\begin{cases}
|\log\lambda|		&	(d=2),\\
|\log\lambda|^{1/2}	&	(d\geq 3).
\end{cases}
$$
This result is valid for any even, convex, not identically zero, external potential $V$ satisfying some mild growth condition (\textit{e.g.}, any polynomial growth is fine). In dimension $1$, on the other hand, the critical behavior does depend on the choice of $V$, see~\cite{HrVe2004}; in this case, it has also been possible to prove exponential decay of covariances.

\subsection{New results}

We consider the measure $\mu^{+,0}_{L;\lambda,\eta}$. Let $\eta>\etac$ and $\lambda>0$. When $\lambda=0$ (\textit{i.e.}, at phase coexistence), the system is in the partial wetting regime, and the interface is expected to be localized. The next theorem shows that this is indeed the case. Moreover, it relates the rate of vanishing of the free energy to the divergence rate of the interface height, as $\eta$ is sent to $\etac$.
\begin{theorem}
\label{thm_UpperBounds}
There exist $T_0>0$, $\bar\eta>\etac$, $\lambda_0>0$, $\alpha_d>0$ and $C_d<\infty$ such that, for any $T>T_0$, $\eta\in(\etac,\bar\eta)$, $\lambda\in(0,\lambda_0)$ and $L\geq 1$,
$$
\mu^{+,0}_{L;\lambda,\eta}\bigl( \varphi_i \geq T|\log f(\lambda,\eta)| \bigr) \leq
C_2 \exp\left(-\alpha_2 T^2 |\log f(\lambda,\eta)|^2/(\log T+|\log f(\lambda,\eta)|)\right),
$$
for $d=2$, and
$$
\mu^{+,0}_{L;\lambda,\eta}\bigl( \varphi_i \geq T\sqrt{|\log f(\lambda,\eta)|} \bigr) \leq
C_d \exp\left(-\alpha_d T^2 |\log f(\lambda,\eta)|\right),
$$
for $d\geq 3$. In particular, there exist $c'_d<\infty$ such that, for all $\eta\in(\etac,\bar\eta)$ and all $L\geq 1$,
$$
\bk{ \varphi_0 }^{+,0}_{L;0,\eta} \leq c'_2 |\log f(0,\eta)|,
$$
when $d=2$, while
$$
\bk{ \varphi_0 }^{+,0}_{L;0,\eta} \leq c'_d \sqrt{|\log f(0,\eta)|},
$$
when $d\geq 3$.
\end{theorem}
\begin{remark}
Since $\bk{ \varphi_0 }^{+,0}_{L;\lambda,\eta}$ is non-increasing in $\eta$ and in $\lambda$ (by FKG inequality), it follows, for example, that
$$
\sup_{\substack{\lambda\geq 0\\L\geq 1}}\bk{ \varphi_0 }^{+,0}_{L;\lambda,\eta} <\infty,
$$
for any $\eta>\etac$ and in any dimension $d\geq 2$.
\end{remark}
We expect that the $|\log f(0,\eta)|$ and $\sqrt{|\log f(0,\eta)|}$ upper bounds are of the correct order. We now state lower bounds of this type.
\begin{theorem}
\label{thm_LowerBounds}
Let $\alpha>1$.
There exists $\bar\eta>\etac$ and $c''_d=c''_d(\alpha)>0$ such that, for all $\eta\in(\etac,\bar\eta)$ and all $L\geq L_0(\eta,\alpha)$,
$$
\bk{ \varphi_0 }^{+,0}_{L;0,\eta} \geq c''_2 |\log f(0,\alpha\eta)|,
$$
when $d=2$, while
$$
\bk{ \varphi_0 }^{+,0}_{L;0,\eta} \geq c''_d \sqrt{|\log f(0,\alpha\eta)|},
$$
when $d\geq 3$.
\end{theorem}
\begin{remark}
Although these results are interesting, it would be more informative to have estimates of the height that are expressed directly in terms of the microscopic parameter $\eta$, and not in terms of the free energy. To do this, one needs to understand the dependence of the latter on $\eta$ close to the wetting transition, a task that seems too hard for the moment. It is however possible to extract a \emph{lower bound} of this type from the proof in~\cite{BoDeZe2000}, which shows that, for $d\geq 3$,
\begin{equation}
\label{eq_feta}
f(0,\eta) \geq c_1(d) e^{-c_2/\eta},
\end{equation}
for some constants $0<c_1,c_2<\infty$. This, combined with the above estimates, implies that
$$
\bk{ \varphi_0 }^{+,0}_{L;0,\eta} \leq c_3 \eta^{-1/2},
$$
for some constant $c_3<\infty$. Observe that if, as we believe, the estimate~\eqref{eq_feta} is of the correct order, then the rate of divergence of the interface height is much faster than the logarithmic divergences seen in the results described in the previous subsection. This would of course be due to the very low density of pinned sites as $\eta$ gets close to $\etac$.
\end{remark}
\begin{remark}
Observe also that the introduction of the parameter $\alpha$ in the lower bound should be irrelevant. Indeed, if the logarithm of the free energy behaves (as indicated by the lower bound of~\cite{BoDeZe2000} for $d\geq 3$) like a polynomial function of $1/\eta$, as $\eta\downarrow \etac$, our upper and lower bounds would actually differ only by a multiplicative constant.
\end{remark}

\subsection{Open problems}
Even though the results presented here substantially improve the description of the wetting transition in these effective models, a number of open problems remain.
\begin{itemize}
\item It would also be desirable to remove the factor $\alpha$ in Theorem~\ref{thm_LowerBounds}. As remarked above, we expect that the latter plays no role, but the verification of this hinges on the next open problem.
\item Obtain information on the behavior of the free energy as a function of $\eta$ close to the wetting transition. This seems too hard at the present time when $d=2$, but there might be some way to prove upper bounds when $d\geq 3$. At a future stage, it would of course be extremely interesting to determine the critical exponent describing the divergence of the height.
\item Prove that the covariances are exponentially decaying with the distance. It is not clear how this should be tackled. The only nonperturbative methods to prove this type of result we are aware of apply only when a suitable graphical representation is available (a random walk representation, for example). However, all the representations available for this model only apply to $2$-point functions, not covariances, and therefore do not seem very helpful.
\item The present work deals only with the partial wetting regime. The situation concerning the complete wetting regime is still not as satisfactory as we would like. In particular, it would be quite desirable to prove that, in the whole complete wetting regime (or at least in the interior of this domain), the interface height diverges like $\log L$ under the measure $\mu^{+,0}_{L;0,\eta}$ in dimension $2$. This is only known to hold far from the critical point. Of course, one expects even more, namely that the fields under $\mu^{+,0}_{L;0,\eta}$ and $\mu^{+,0}_{L;0,0}$ should be very close.
\end{itemize}

\begin{acknowledgements}
The author would like to express his warmest thanks to an anonymous referee for pointing out numerous mistakes and imprecisions in the submitted version of this work, and for suggesting many improvements to the exposition.
Discussions with Erwin Bolthausen, Jean-Dominique Deuschel and Dima Ioffe are also gratefully acknowledged.
\end{acknowledgements}
\section{Proofs}
The main tools used in the proofs below are FKG inequality and the chessboard estimate. Both hold for the measures considered here, because they do for the Gaussian measure, and are insensitive to perturbation of the form $\prod_i e^{U(\varphi_i)}$ (after suitably smoothing our potential -- in particular the positivity constraint and the pinning potential -- and taking weak limits). We refer, \textit{e.g.}, to Appendix~B of~\cite{Gi2001} for additional information and references on the validity and use of FKG inequality in the context of gradient fields, and to~\cite{Bi2006} for a nice review on reflection positivity (and, in particular, the chessboard estimate). Let us also emphasize, to make the following arguments clearer, that in all the applications of the chessboard estimates in the present work, we are using reflection through planes between lattice sites.

\medskip
In order to use reflection positivity, we shall need to work with periodic boundary condition.
Let us quickly recall the corresponding definitions. We denote by $\bbT^d_L$ the torus $\Zd/(L\Zd)$. Configurations are then given by $\varphi\in\bbR^{\bbT^d_L}$. The measures $\mu^{\per}_{L;\lambda,\eta}$ and $\mu^{+,\per}_{L;\lambda,\eta}$ are defined precisely as before,
\begin{align*}
\mu^{\per}_{L;\lambda,\eta} (\dd\varphi) &= \frac1{\PF^{\per}_{L;\lambda,\eta}} \exp \bigl[-\frac1{8d}\sum_{\substack{i,j\in\bbT^d_L\\i\sim j}} (\varphi_i-\varphi_j)^2 - \lambda \sum_{i\in\bbT^d_L} V(\varphi_i) \bigr]\\
&\hspace{5cm}\times \prod_{i\in\bbT^d_L} \bigl( \dd\varphi_i + \eta\delta_0(\dd\varphi_i) \bigr),\\
\mu^{+,\per}_{L;\lambda,\eta} (\,\cdot\,) &= \mu^{\per}_{L;\lambda,\eta} (\,\cdot\given \Omega_+),
\end{align*} 
reinterpreting $i\sim j$ to mean that $i$ and $j$ are neighboring vertices on $\bbT^d_L$. Notice that for these measures to be well-defined, it is necessary that $\lambda>0$.

We denote by $\PF^{\per}_{L;\lambda,\eta}$ and $\PF^{+,\per}_{L;\lambda,\eta}$ the corresponding partition functions, and by
$$
f^\per_L(\lambda,\eta) = L^{-d} \log(\PF^{+,\per}_{L;\lambda,\eta}/\PF^{+,\per}_{L;\lambda,0})
$$
the corresponding free energy. In the thermodynamic limit, this free energy and the one defined with $0$-boundary condition agree.
\begin{lemma}
\label{lem_FreeEnergy}
For all $\lambda\geq 0,\eta\geq 0$, the limit $f(\lambda,\eta) = \lim_{L\to\infty}f^0_L(\lambda,\eta)$ exists and is convex and increasing in $\eta$ and $\lambda$. Moreover, for all $\lambda>0,\eta\geq 0$, the limit $\lim_{L\to\infty}f^\per_L(\lambda,\eta)$ also exists and coincides with $f(\lambda,\eta)$.
\end{lemma}
\begin{proof}
The existence of $\lim_{L\to\infty}f^0_L(\lambda,\eta)$ follows, for example, by FKG inequality and completely standard superadditivity arguments. Its monotonicity in $\lambda$ is also immediate
$$
\frac{\partial}{\partial\lambda} f^0_L(\lambda,\eta) = L^{-d} \sum_{i\in\Lambda_L} \bigl\{ \bk{V(\varphi_i)}^{+,0}_{L;\lambda,0} - \bk{V(\varphi_i)}^{+,0}_{L;\lambda,\eta} \bigr\} \geq 0,
$$
by FKG inequality, since $V$ is increasing on $\bbR^+$. To check the monotonicity in $\eta$, it suffices to observe that $\eta\partial/\partial\eta f^0_L(\lambda,\eta)$ is simply the density of pinned sites, and thus positive. Convexity follows similarly, since second derivatives yield variances.

Let us denote by $\PF[A]$ the restriction of the partition function $\PF$ to configurations satisfying the condition $A$. With this notation, we have,
\begin{equation}
\label{eq_almostbouded}
\tfrac12\PF^{+,\per}_{L;\lambda,\eta} \leq \PF^{+,\per}_{L;\lambda,\eta}\bigl[\varphi_i \leq V^{-1}((2d/\lambda)\log L), \forall i\in\bbT^d_L\bigr] \leq \PF^{+,\per}_{L;\lambda,\eta},
\end{equation}
for all $L$ large enough.
Indeed,
$$
\frac{\PF^{+,\per}_{L;\lambda,\eta}\bigl[\varphi_i \leq V^{-1}(\frac{2d}{\lambda}\log L), \forall i\in\bbT^d_L\bigr]} {\PF^{+,\per}_{L;\lambda,\eta}} = \mu^{+,\per}_{L;\lambda,\eta}(\varphi_i \leq V^{-1}(\frac{2d}{\lambda}\log L), \forall i\in\bbT^d_L),
$$
which proves the second inequality, and, by FKG inequality,
$$
\mu^{+,\per}_{L;\lambda,\eta}(\varphi_i \leq V^{-1}(\frac{2d}{\lambda}\log L), \forall i\in\bbT^d_L)) \geq \prod_{i\in\bbT^d_L} \mu^{+,\per}_{L;\lambda,\eta}(\varphi_i \leq V^{-1}(\frac{2d}{\lambda}\log L)).
$$
But another application of FKG, and the chessboard estimate yield
\begin{align*}
\mu^{+,\per}_{L;\lambda,\eta}(\varphi_i > V^{-1}(\frac{2d}{\lambda}\log L))
&\leq
\mu^{+,\per}_{L;\lambda,0}(\varphi_i > V^{-1}(\frac{2d}{\lambda}\log L))\\
&\leq
\mu^{+,\per}_{L;\lambda,0}(\varphi_j > V^{-1}(\frac{2d}{\lambda}\log L), \forall j\in\bbT^d_L)^{1/|\bbT^d_L|}\\
&\leq
\exp\bigl( -\lambda V(V^{-1}(\frac{2d}{\lambda}\log L)) \bigr)\\
&=
L^{-2d}.
\end{align*}
Therefore,
$$
\mu^{+,\per}_{L;\lambda,\eta}(\varphi_i \leq V^{-1}((2d/\lambda)\log L), \forall i\in\bbT^d_L)) \geq e^{-2 L^{-d}} \geq \frac12,
$$
for all $L$ large enough. This proves~\eqref{eq_almostbouded}. Notice now that the same also holds for $0$-boundary condition. Indeed, the upper bound is again trivial, and for the lower bound, we can use FKG inequality to get
$$
\mu^{+,0}_{L;\lambda,\eta}(\varphi_i \leq V^{-1}(\frac{2d}{\lambda}\log L), \forall i\in\bbT^d_L))
\geq \mu^{+,\per}_{L+1;\lambda,\eta}(\varphi_i \leq V^{-1}(\frac{2d}{\lambda}\log L), \forall i\in\bbT^d_L)).
$$
It is thus sufficient to compare the partition function with periodic and $0$-boundary conditions, under the constraint that all spins satisfy $\varphi_i<V^{-1}(\frac{2d}{\lambda}\log L)$. However, for a configuration $\varphi$ on the torus satisfying this constraint, the change in energy resulting from setting one height to $0$ is bounded above by $\frac12(V^{-1}(\frac{2d}{\lambda}\log L))^2$ and bounded below by
$$
-\frac12(V^{-1}(\frac{2d}{\lambda}\log L))^2- 2d\log L.
$$
Since $|\partial\Lambda_L|=O(L^{d-1})$, this shows that, for all fixed $\lambda>0$,
$$
L^{-d}\, \left| \log
\frac{\PF^{+,0}_{L;\lambda,\eta}\bigl[\varphi_i \leq V^{-1}(\frac{2d}{\lambda}\log L), \forall i\in\bbT^d_L)\bigr]}{\PF^{+,\per}_{L+1;\lambda,\eta}\bigl[\varphi_i \leq V^{-1}(\frac{2d}{\lambda}\log L), \forall i\in\bbT^d_L\bigr]} \right| = o_L(1),
$$
implying that the limiting free energies coincide. By the above considerations, this is also true for the unrestricted partition functions and free energies. This concludes the proof of Lemma~\ref{lem_FreeEnergy}.
\end{proof}

\subsection{Upper bound on the height: proof of Theorem~\ref{thm_UpperBounds}}
In this section, we prove the upper tail estimate for the height at the origin, and the resulting upper bound on the height of the interface. 

In the proof, it will be convenient to assume from the start that the free energy is small enough; more precisely, we shall always assume that $\bar\eta$ and $\lambda_0$ are chosen in such a way that $f(\lambda_0,\bar\eta)\leq e^{-1}$, which ensures that $|\log f(\lambda,\eta)|\geq 1$ and $f(\lambda,\eta)^{-1/d}\geq 1$, for all $\eta\in(\etac,\bar\eta)$ and all $\lambda\in(0,\lambda_0)$.

Let us first observe that FKG inequality implies that $\lim_{L\to\infty}\bk{\varphi_0}^{+,0}_{L;\lambda,\eta}$ exists in $\bbR\cup\{+\infty\}$. It is therefore sufficient for us to restrict our attention to boxes of size $L+1=2^N$, $N\geq 0$, when $d=2$, and $L=2^N$, $N\geq 0$, when $d\geq 3$. This will be useful to ensure that the sizes of the blocks used when applying the chessboard estimate divide the size of the torus.

\subsubsection{The two-dimensional case}
Let us fix $\lambda>0$ and $\eta>\etac$ as above.
Expanding over pinned sites (see~\cite{DeVe2000,IoVe2000,BoVe2001}, for example), we have
\begin{multline*}
\mu^{+,0}_{L;\lambda,\eta}\bigl(\varphi_0 \geq T|\log f(\lambda,\eta)|\bigr)\\
=
\sum_{k\geq 1} \sum_{\substack{A\cap B_{k-1} = \emptyset\\A\cap B_{k} \neq \emptyset}} \pinned^{+,0}_{L;\lambda,\eta}(A)\, \mu^{+,0}_{L;\lambda,0}\bigl(\varphi_0 \geq T|\log f(\lambda,\eta)| \bgiven \varphi_i=0, \forall i\in A \bigr),
\end{multline*}
where we used the notation $B_k=\setof{i\in\bbT^d_L}{\normINF{i}\leq k}$, $k\geq 0$.

By FKG and Lemma~\ref{lem_singlePinned} below (provided that $T$ is large enough), we can find $c_1$ such that
\begin{multline*}
\mu^{+,0}_{L;\lambda,0}\bigl(\varphi_0 \geq T|\log f(\lambda,\eta)| \bgiven \varphi_i=0, \forall i\in A \bigr)\\
\leq
\sup_{\normINF{i}=k} \mu^{+,0}_{L;0,0}\bigl(\varphi_0 \geq T|\log f(\lambda,\eta)| \bgiven \varphi_i=0 \bigr)\\
\leq
\exp\bigl(-c_1 T^2|\log f(\lambda,\eta)|^2 / \log k\bigr),
\end{multline*}
uniformly in $A$ such that $A\cap B_{k-1} = \emptyset$ and $A\cap B_{k} \neq \emptyset$. This implies that
$$
\mu^{+,0}_{L;\lambda,\eta} \bigl( \varphi_0\geq T|\log f(\lambda,\eta)| \bigr)
\leq
\sum_{k\geq 1} \exp\bigl(-c_1 T^2|\log f(\lambda,\eta)|^2 / \log k\bigr) \sum_{A\cap B_{k-1} = \emptyset} \pinned^{+,0}_{L;\lambda,\eta}(A).
$$
But, for all $k\geq 1$,
$$
\sum_{A\cap B_{k-1} = \emptyset} \pinned^{+,0}_{L;\lambda,\eta}(A)
=
\frac{\PF^{+,0}_{L;\lambda,\eta}(B_{k-1}^\comp)}{\PF^{+,0}_{L;\lambda,\eta}},
$$
where $\PF^{+,0}_{L;\lambda,\eta}(B_{k-1}^\comp)$ is defined as in~\eqref{equ_GibbsMeasure} but with the pinning potential acting only on $B_{k-1}^\comp$.

To estimate this last ratio, we would like to follow the idea in~\cite{BoBr2001} and use reflexion positivity of the Gibbs measure (which holds, since we are considering a nearest-neighbor gradient field, with on-site potentials). Of course, the pinning potential is a bit singular, and makes the application of this inequality awkward, so we first replace it by its more regular approximation~\eqref{eq_SWpinning}. The only remaining obstacle now is that we have $0$-boundary condition instead of periodic boundary conditions. To remove this problem, we use once more FKG inequality to obtain, for any $k\geq 1$,
\begin{align*}
\frac{\PF^{+,0,(\epsilon)}_{L;\lambda,\eta}} {\PF^{+,0,(\epsilon)}_{L;\lambda,\eta}(B_{k-1}^\comp)}
&=
\bk{\prod_{i\in B_{k-1}} e^{-U^{(\epsilon)}_{2\eta}(\varphi_i)}}_{L;\lambda,\eta}^{+,0,(\epsilon)}(B_{k-1}^\comp)\\
&\geq
\bk{\prod_{i\in B_{k-1}} e^{-U^{(\epsilon)}_{2\eta}(\varphi_i)}}_{L+1;\lambda,\eta}^{+,\per,(\epsilon)}(B_{k-1}^\comp)
= \frac{\PF^{+,\per,(\epsilon)}_{L+1;\lambda,\eta}}{\PF^{+,\per,(\epsilon)}_{L+1;\lambda,\eta}(B_{k-1}^\comp)}.
\end{align*}
Let $n\in\bbN$ be such that $2^n\leq k < 2^{n+1}$, and set $\bar{R}=2^{n-1}$. Applying the chessboard estimate starting with the block $\{-\bar{R}+1,\ldots,\bar{R}\}^2$, we get
\begin{align*}
\frac{\PF^{+,\per,(\epsilon)}_{L+1;\lambda,\eta}(B_{k-1}^\comp)} {\PF^{+,\per,(\epsilon)}_{L+1;\lambda,\eta}}
&=
\bk{\prod_{i\in B_{k-1}} e^{U^{(\epsilon)}_{2\eta}(\varphi_i)}}_{L+1;\lambda,\eta}^{+,\per,(\epsilon)}\\
&\leq
\bk{\prod_{i\in B_{\bar{R}}} e^{U^{(\epsilon)}_{2\eta}(\varphi_i)}}_{L+1;\lambda,\eta}^{+,\per,(\epsilon)}\\
&\leq
\bigl\{ \bk{\prod_{i\in\bbT^2_{L+1}} e^{U^{(\epsilon)}_{2\eta}(\varphi_i)}}_{L+1;\lambda,\eta}^{+,\per,(\epsilon)} \bigr\}^{|B_{\bar R}|/|\bbT^2_{L+1}|}\\
&\leq
\bigl\{ \frac{\PF^{+,\per,(\epsilon)}_{L+1;\lambda,0}} {\PF^{+,\per,(\epsilon)}_{L+1;\lambda,\eta}} \bigr\}^{|B_{\bar R}|/|\bbT^2_{L+1}|}\\
\end{align*}
We now obtain the result for our original ratio by taking the limit $\epsilon\downarrow 0$, which yields
\begin{align*}
\frac{\PF^{+,0}_{L;\lambda,\eta}(B_{k-1}^\comp)}{\PF^{+,0}_{L;\lambda,\eta}}
&\leq
\bigl\{ \frac{\PF^{+,\per}_{L+1;\lambda,0}} {\PF^{+,\per}_{L+1;\lambda,\eta}} \bigr\}^{|B_{\bar R}|/|\bbT^2_{L+1}|}\\
&=
e^{-|B_{\bar R}|\, f_{L+1}^{\per}(\lambda,\eta)}\\
&\leq
e^{-\tfrac14 |B_{k-1}|\, f_{L+1}^{\per}(\lambda,\eta)}.
\end{align*}
Collecting these estimates, and using Lemma~\ref{lem_FreeEnergy}, we obtain finally that
$$
\sum_{A\cap B_{k-1} = \emptyset} \pinned^{+,0}_{L;\lambda,\eta}(A) \leq e^{-\tfrac14 |B_{k-1}|\, (f(\lambda,\eta)-o_L(1))} \leq e^{-c_2 f(\lambda,\eta) k^2},
$$
for all $L>L_0(\lambda,\eta)$ large enough. Consequently, setting
$$
\bar{k}=\frac{T|\log f(\lambda,\eta)|}{\sqrt{f(\lambda,\eta)(\log T + |\log f(\lambda,\eta)|)}},
$$
we obtain the bound, valid for all $L>L_0(\lambda,\eta)$,
\begin{align*}
\mu^{+,0}_{L;\lambda,\eta} \bigl( \varphi_0\geq T|\log f(\lambda,\eta)| \bigr)
&\leq
\sum_{k\geq 1} e^{-c_1 T^2|\log f(\lambda,\eta)|^2 / \log k}e^{- c_2 f(\lambda,\eta) k^2}\\
&\leq
\sum_{k=1}^{\bar k} e^{-c_1 T^2|\log f(\lambda,\eta)|^2 / \log k} + \sum_{k>\bar k} e^{- c_2 f(\lambda,\eta) k^2}\\
&\leq \bar{k} e^{-c_1 T^2|\log f(\lambda,\eta)|^2 / \log \bar k} + e^{-c_2 f(\lambda,\eta) \bar k^2} \sum_{k=1}^\infty e^{-c_2 f(\lambda,\eta) k^2}\\
&\leq
c_3 e^{-c_4 T^2|\log f(\lambda,\eta)|^2 / (\log T + |\log f(\lambda,\eta)|)}.
\end{align*}
However, this probability is increasing in $L$ (by FKG), so that the above bound actually holds for all $L\geq 1$. We can now also easily deduce the stated upper bound on the mean height: For all $L\geq 1$,
\begin{align*}
\bk{\varphi_0}^{+,0}_{L;\lambda,\eta}
&\leq T_0 |\log f(\lambda,\eta)|
+ \int_{T_0|\log f(\lambda,\eta)|}^\infty \mu^{+,0}_{L;\lambda,\eta} \bigl( \varphi_0\geq t \bigr) \dd t\\
&\leq
c_5 |\log f(\lambda,\eta)|.
\end{align*}
In particular, taking the limit as $\lambda\downarrow 0$, we get
$$
\bk{\varphi_0}^{+,0}_{L;0,\eta}
\leq
c_5 |\log f(0,\eta)|.
$$

All statements have been proved now, except for the following lemma, which was used in the above argument.
\begin{lemma}
\label{lem_singlePinned}
There exists $T_0$ and $c>0$ such that, for all $T>T_0$ and $L> R\geq 1$,
$$
\sup_{\normINF{i}= R}\mu^{+,0}_{L;0,0}\bigl(\varphi_0 \geq T \log R \bgiven \varphi_i=0 \bigr) \leq R^{-cT^2}.
$$
\end{lemma}
\begin{proof}
The proof is a variant of the one given in~\cite[Lemma~4.4]{Gi2001}. Let $i$ be such that $\normINF{i}=R$. We first use FKG inequality to center the box around $i$ (at the cost of replaing $L$ by $2L$),
$$
\mu^{+,0}_{L;0,0}\bigl(\varphi_0 \geq T \log R \bgiven \varphi_i=0  \bigr) \leq \mu^{+,0}_{2L;0,0}\bigl(\varphi_{-i} \geq T \log R \bgiven \varphi_0=0  \bigr).
$$
Using again FKG inequality, we can also deduce that the latter probability only increases if we replace the $0$-boundary condition with $\alpha\log L$-boundary condition, $\alpha>0$,
$$
\mu^{+,0}_{2L;0,0}\bigl(\varphi_{-i} \geq T \log R \bgiven \varphi_0=0  \bigr) \leq
\mu^{+,\alpha\log L}_{2L;0,0}\bigl(\varphi_{-i} \geq T \log R \bgiven \varphi_0=0 \bigr).
$$
Choosing $\alpha$ large enough, one can guarantee that $\mu^{+,\alpha\log L}_{2L;0,0}\bigl(\Omega^+ \bgiven \varphi_0=0 \bigr) > \tfrac12$, see~\cite[Lemma~4.4]{Gi2001}. Consequently, we can remove the positivity constraint,
$$
\mu^{+,\alpha\log L}_{2L;0,0}\bigl(\varphi_{-i} \geq T \log R \bgiven \varphi_0=0 \bigr) \leq
2 \mu^{\alpha\log L}_{2L;0,0}\bigl(\varphi_{-i} \geq T \log R \bgiven \varphi_0=0 \bigr).
$$
Now, under $\mu^{\alpha\log L}_{2L;0,0}\bigl(\,\cdot \bgiven \varphi_0=0 \bigr)$, $\varphi_{-i}$ is Gaussian and the random walk representation implies that its mean is at most $c'\alpha\log R$, while its variance is at least $c''\log R$. Therefore, provided that $T$ is large enough,
$$
\mu^{\alpha\log L}_{2L;0,0}\bigl(\varphi_{-i} \geq T \log R \bgiven \varphi_0=0 \bigr) \leq \exp\bigl( -c T^2 \log R \bigr),
$$
which proves the claim.
\end{proof}
\subsubsection{Dimension $3$ and higher}
The argument in dimensions $3$ and higher is unfortunately more involved, due to the fact that pinning a single point does not localize the interface anymore. One way to solve this problem would be to use reflection positivity to deduce that the distribution of pinned sites dominates (in a suitable sense) some Bernoulli percolation-type process, and then study the entropic repulsion problem for an interface with a such a random distribution of pinned sites. This certainly looks feasible, but we decided to try a different path relying more on reflection positivity, but in our opinion technically simpler.

We start in a way similar to what we did in dimension $2$. We fix $\lambda>0$ and $\eta>\etac$ so that the free energy is sufficiently small. We write $R(k) = \lfloor k\,f(\lambda,\eta)^{-1/d} \rfloor$. For $k\geq 0$, let $\calE_1(k) = \{\varphi_j>0,\, \forall j\in B_{R(k)}\}$ and $\calE_2(k) = \{\exists i\in B_{R(k+1)}\setminus B_{R(k)},\, \varphi_i=0\}$. It follows from Cauchy-Schwarz inequality that
\begin{align*}
\mu^{+,\per}_{L;\lambda,\eta}(\varphi_0 \geq T&\sqrt{|\log f|})\\
&=
\sum_{k\geq 0} \mu^{+,\per}_{L;\lambda,\eta}(\varphi_0 \geq T\sqrt{|\log f|}, \calE_1(k), \calE_2(k))\\
&\leq
\sum_{k\geq 0} \Bigl[ \mu^{+,\per}_{L;\lambda,\eta}(\varphi_0 \geq T\sqrt{|\log f|}, \calE_2(k))\, \mu^{+,\per}_{L;\lambda,\eta}(\calE_1(k)) \Bigr]^{1/2}.
\end{align*}
First, observe that applying the chessboard estimate similarly as was done in the two-dimensional case, we obtain
\begin{align*}
\mu^{+,\per}_{L;\lambda,\eta}(\calE_1(k))
&\leq
\left[\mu^{+,\per}_{L;\lambda,\eta}(\varphi_j>0,\, \forall j\in \bbT^d_L)\right]^{c_1 \left(\frac {R(k)}L\right)^d}\\
&\leq
\left[\frac{\PF^{+,\per}_{L;\lambda,0}}{\PF^{+,\per}_{L;\lambda,\eta}}\right]^{c_1 \left(\frac {R(k)}L\right)^d}\\
&\leq
\exp\bigl( -c_2 (f(\lambda,\eta)-o_L(1)) {R(k)}^d \bigr)\\
&\leq
\exp\bigl( -c_3\, k^d \bigr),
\end{align*}
for all $L$ large enough.

It remains to estimate $\mu^{+,\per}_{L;\lambda,\eta}(\varphi_0 \geq T\sqrt{|\log f|}, \calE_2(k))$. To lighten notations, let us write $h=T\sqrt{|\log f|}$. Let $n\in\bbN$ be such that $2^n\geq 2R(k)+1 > 2^{n-1}$, and set $\bar{R}=2^{n}$. Let us fix $i\in B_{R(k+1)}\setminus B_{R(k)}$. Applying once more the chessboard estimate, starting this time with the block $i+\{-\bar{R}+1,\ldots,\bar{R}\}^d$, yields
\begin{equation*}
\mu^{+,\per}_{L;\lambda,\eta}(\varphi_0 \geq h, \varphi_i=0)
\leq
\Bigl[\mu^{+,\per}_{L;\lambda,\eta}\bigl(\varphi_j \geq h,\forall j\in\calR(0), \varphi_k=0, \forall k\in\calR(i)\bigr)\Bigr]^{|B_{\bar{R}}|/|\bbT^d_L|},
\end{equation*}
where, for $j$ a site of the original block, $\calR(j)$ is the set of sites obtained after applying all the reflections to $j$.

We can now get rid of the mass and the pinning potential. Using FKG inequality,
\begin{align*}
\mu^{+,\per}_{L;\lambda,\eta}\bigl(\varphi_j \geq h,\forall j\in\calR(0)&, \varphi_k=0, \forall k\in\calR(i)\bigr)\\
&\leq
\mu^{+,\per}_{L;\lambda,\eta}\bigl(\varphi_j \geq h,\forall j\in\calR(0) \given \varphi_k=0, \forall k\in\calR(i)\bigr)\\
&\leq
\mu^{+,\per}_{L;0,0}\bigl(\varphi_j \geq h,\forall j\in\calR(0) \given \varphi_k=0, \forall k\in\calR(i)\bigr).
\end{align*}
For technical reasons, we shall need the distance between the points of $\calR(i)$ to be larger than $R_0$, for some constant $R_0=R_0(d)$ to be fixed later. In order to do that, let us introduce $\Delta(d)=\min\setof{n\geq 1}{2n\bar{R}\geq R_0}$ and let us denote by $\calR_\Delta(i)$ the subset of $\calR(i)$ obtained by replacing $\bar{R}$ by $\Delta\bar{R}$ in the construction. Of course, the distance between points in $\calR_\Delta(i)$ is at least $R_0$. By FKG inequality,
\begin{multline*}
\mu^{+,\per}_{L;0,0}\bigl(\varphi_j \geq h,\forall j\in\calR(0) \bgiven \varphi_k=0, \forall k\in\calR(i)\bigr)\\
\leq
\mu^{+,\per}_{L;0,0}\bigl(\varphi_j \geq h,\forall j\in\calR(0) \bgiven \varphi_k=0, \forall k\in\calR_\Delta(i)\bigr).
\end{multline*}
We want now to deal with the positivity constraint.
\begin{multline*}
\mu^{+,\per}_{L;0,0}\bigl(\varphi_j \geq h,\forall j\in\calR(0) \bgiven \varphi_k=0, \forall k\in\calR_\Delta(i)\bigr)\\
\leq
\frac{\mu^{\per}_{L;0,0}\bigl(\varphi_j \geq h,\forall j\in\calR(0) \bgiven \varphi_k=0, \forall k\in\calR_\Delta(i)\bigr)}{\mu^{\per}_{L;0,0}\bigl(\Omega_+ \bgiven \varphi_k=0, \forall k\in\calR_\Delta(i)\bigr)}.
\end{multline*}
The denominator in the last expression has already been bounded below in~\cite{BoDeZe2000} in the case of $0$-boundary condition (it is here that we need $2\Delta\bar{R}\geq R_0$). Using FKG inequality to change correspondingly the boundary condition, this yields
\begin{align*}
\mu^{\per}_{L;0,0}\bigl(\Omega_+ \given \varphi_k=0, \forall k\in\calR_\Delta(i)\bigr)
&\geq
\mu^{0}_{L-1;0,0}\bigl(\Omega_+ \given \varphi_k=0, \forall k\in\calR_\Delta(i)\bigr)\\
&\geq
\exp(-c_4\, (L/(\Delta\bar{R}))^d \log(\Delta\bar{R}))\\
&\geq
\exp(-c_5\, (L/R(k))^d \log R(k)).
\end{align*}
Obviously,
\begin{multline*}
\mu^{\per}_{L;0,0}\bigl(\varphi_j \geq h,\forall j\in\calR(0) \bgiven \varphi_k=0, \forall k\in\calR_\Delta(i)\bigr)\\
\leq \mu^{\per}_{L;0,0}\bigl(\sum_{j\in\calR(0)}\varphi_j \geq h|\calR(0)| \bgiven \varphi_k=0, \forall k\in\calR_\Delta(i)\bigr).
\end{multline*}
Under $\mu^{\per}_{L;0,0}\bigl(\,\cdot \bgiven \varphi_k=0, \forall k\in\calR_\Delta(i)\bigr)$, $\sum_{j\in\calR(0)}\varphi_j$ is a Gaussian random variable of mean $0$, and variance
$$
\sum_{j,k\in\calR(0)} \sum_{n\geq 0} \bbP_j(X_n=k, \tau_{\calR_\Delta(i)}>n),
$$
where $\bbP_j$ is the law of the simple random walk $(X_n)_{n\geq 0}$ on $\bbT^d_L$ with $X_0=j$, and $\tau_{\calR(i)}=\min\setof{n\geq 0}{X_n\in\calR_\Delta(i)}$. Lemma~\ref{lem_twogrids} below implies that this variance is bounded above by $c_6 |\calR(0)|$,
and thus
\begin{equation*}
\mu^{\per}_{L;0,0}\bigl(\varphi_j \geq h,\forall j\in\calR(0) \given \varphi_k=0, \forall k\in\calR_\Delta(i)\bigr)
\leq \exp(-c_7\, (L/R(k))^d h^2),
\end{equation*}
Putting all together, we have shown that
$$
\mu^{+,\per}_{L;\lambda,\eta}(\varphi_0 \geq h, \varphi_i=0) \leq \exp(-c_8\, T^2|\log f(\lambda,\eta)|),
$$
provided that $\log k \leq c_8 T^2|\log f(\lambda,\eta)|\equiv \log k_{\scriptscriptstyle\rm max}$, with $c_8=c_7/4c_5$, and $T\geq \sqrt{2c_5/dc_7}\equiv T_1$. Thus,
\begin{align*}
\mu^{+,\per}_{L;\lambda,\eta}(\varphi_0 \geq T\sqrt{|\log f(\lambda,\eta)|}, \calE_2(k))
&\leq
 \exp(-c_8\, T^2|\log f(\lambda,\eta)|),
\end{align*}
for all $k<k_{\scriptscriptstyle\rm max}$.

We have proved that, for all $L\geq L_0(\eta,\lambda)$,
\begin{align*}
\mu^{+,\per}_{L;\lambda,\eta}(\varphi_0 \geq T\sqrt{|\log f(\lambda,\eta)|})
&\leq
 \exp(-c_8\, T^2|\log f(\lambda,\eta)|) \sum_{k=1}^{k_{\scriptscriptstyle\rm max}} e^{- c_3 k^d}
+
\sum_{k>k_{\scriptscriptstyle\rm max}} e^{- c_3 k^d}\\
&\leq
c_{11}  \exp(-c_8\, T^2|\log f(\lambda,\eta)|),
\end{align*}
provided that $T\geq T_1$.

Proceeding as in the two-dimensional case, one finally obtains that, for all $L\geq L_0(\eta,\lambda)$,
$$
\bk{\varphi_0}^{+,\per}_{L;\lambda,\eta}
\leq
c_{12} \sqrt{|\log f(\lambda,\eta)|}.
$$
The results stated for $0$-boundary condition follows from FKG inequality (first to change the boundary condition, and then to argue as in the two-dimensional case).

To complete the proof, it only remains to prove the following lemma.
\begin{lemma}
\label{lem_twogrids}
For $\calR(0)$ and $\calR_\Delta(i)$ defined as above,
$$
\sum_{j,k\in\calR(0)} \sum_{n\geq 0} \bbP_j(X_n=k, \tau_{\calR_\Delta(i)}>n) \leq c_6|\calR(0)|.
$$
\end{lemma}
\begin{proof}
Of course,
$$
\sum_{j,k\in\calR(0)} \sum_{n\geq 0} \bbP_j(X_n=k, \tau_{\calR_\Delta(i)}>n)
=
\sum_{j\in\calR(0)} \sum_{n\geq 0} \bbP_j(X_n\in\calR(0), \tau_{\calR_\Delta(i)}>n).
$$
\smallskip
We start by periodizing the sets $\calR_\Delta(i)$ and $\calR(0)$.

$\calR_\Delta(i)$ contains exactly one image of $i$ in each of the blocks used during its construction; it can thus be partitioned into $2^d$ disjoint periodic subsets of equal sizes. We denote by $S_i$ the subset containing $i$. If we only kill the random walk once it enters $S_i$, then the sum we want to control is only made bigger.

Similarly, the set $\calR(0)$ contains a fixed number of sites in each of these blocks (their number depending on the value of $\Delta$). As before, we can decompose $\calR_0$ as a finite union of disjoint periodic arrays of sites, with the same period as $S_i$. We can of course restrict our attention to one of these subsets only, because if we show separately for each of these subsets that the average number of times they are visited by the walk, before it dies, is bounded above uniformly in everything but the dimension, then the same will be true for their union. Let us therefore consider one of these subsets, which we call $S_0$.

\medskip
To prove the lemma, we are going to show that after each visit of the random walk to the set $S_0$, the walk has a positive probability, depending only on the dimension of the lattice, of hitting $S_i$ before reentering $S_0$. This will show that the number of visits to $\calR(0)$ before entering $S_i$ is stochastically dominated by a geometric random variable of positive parameter (uniformly in everything, but the dimension), which immediately implies the claim.

Now, the periodicity of $R_0$ and $R_i$ allows us to reinterpret the problem as being on a torus of size the common period of these two sets, with two distinguished sites, $s_0$ and $s_i$, coming from $S_0$ and $S_i$, respectively. In these terms, the problem can be reformulated as follows: prove that starting from $s_0$, the random walk has a positive probability of hitting $s_1$ before returning to $s_0$, and this probability is uniform in the size of the torus.

But this is easy. Indeed, by symmetry, at least one half of the sites of the torus satisfy
$$
\bbP_x (\tau_{s_0} > \tau_{s_i}) \geq 1/2,
$$
where $\tau_y=\min\setof{n\geq 1}{X_n=y}$. Let us call this set $\calG$. It is therefore possible to find $r$ such that at least half of the sites in
$$
\setof{x}{\normINF{x} = r}
$$
belong to $\calG$ (otherwise it would be impossible for $\calG$ to contain at least half of the sites of the torus). But the probability that the random walk starting at $0$ exits the box of radius $r$ at one of the sites belonging to $\calG$ before returning to $0$, is bounded away from $0$, uniformly in everything but the dimension (since (i) the random walk is transient, and (ii) the probabilities that the random walk exits the box through any given site are comparable (see, \textit{e.g.}, Lemma~1.7.4 in~\cite{La1991})). The conclusion follows, for once the walk has reached a site of $\calG$, it hits $s_i$ first with probability at least $1/2$.
\end{proof}
\subsection{Lower bound on the height: proof of Theorem~\ref{thm_LowerBounds}}
Let $\alpha>1$ and set $R=\bigl\lfloor \bigl((K f(0,\alpha\eta)\bigr)^{-1/d} \bigr\rfloor$, where $K$ will be chosen (large enough) later (depending on $\alpha$). Let us define the subset $\Lambda_L^R$ of $\Lambda_L$ by the requirement that $\Lambda_L=\bigcup_{i\in\Lambda_L^R} B_R(i)$, where the boxes $B_R(i)$ (cubes of radius $R$ centered at $i$) are disjoint (we assume, without loss of generality, that $L$ is a suitable multiple of $R$). Let $\calB_R^i$ be the event that there are no pinned sites at distance less than $R$ from $i$, $\calB_R^i=\{\varphi_j>0,\, \forall j\in B_R(i)\}$. We also denote by $\calA$ the (random) set of pinned sites. We first observe that
$$
\mu^{+,0}_{2L;0,\eta}(\calB_R^0) \geq \mu^{+,0}_{L;0,\eta}(\calB_R^i),
$$
for all $i\in\Lambda_L^R$, thanks to FKG inequality. This implies that
\begin{align*}
\mu^{+,0}_{2L;0,\eta}(\calB_R^0) &\geq \bk{|\Lambda_L^R|^{-1} \sum_{i\in\Lambda_L^R} \mathbf{1}_{\calB_R^i}}^{+,0}_{L;0,\eta}\\ &\geq \bk{|\Lambda_L^R|^{-1} \sum_{i\in\Lambda_L^R} \mathbf{1}_{\calB_R^i} \given \calE}^{+,0}_{L;0,\eta} \; \mu^{+,0}_{L;0,\eta}(\calE),
\end{align*}
where $\calE=\{\frac{|\calA|}{|\Lambda_L|} \leq \frac2{\log\alpha}f(0,\alpha\eta)\}$.

Let us now consider the average density of pinned sites,
$$
\rho_L(\eta)=|\Lambda_L|^{-1}\bk{|\calA|}_{L;0,\eta}^{+,0}.
$$
This density can easily be bounded above (just use the standard integration-differentiation trick, or see~\cite{CaVe2000}):
$$
f_{L}^{+,0}(0,\alpha\eta) = \int_0^{\alpha\eta} \frac1{t} \rho_L (t) \,\dd t \geq \int_\eta^{\alpha\eta} \frac1{t} \rho_L (t) \,\dd t \geq  \log(\alpha)\, \rho_L(\eta),
$$
where we used the fact that $\rho_L(\eta)$ is a nonnegative and nondecreasing function of $\eta$. Since, for all $L\geq L_0(\eta,\alpha)$, $f(0,\alpha\eta) \geq \tfrac12 f_{L}^{+,0}(0,\alpha\eta)$, Markov inequality implies that
$$
\mu^{+,0}_{L;0,\eta}(\calE^\comp) \leq \mu^{+,0}_{L;0,\eta}\bigl(|\calA| > 2^{-2d-2}K\log\alpha \rho_L(0,\eta)\,|\Lambda_L|\bigr) \leq \frac12,
$$
as soon as $K>2^{2d+3}/\log\alpha$ and $L\geq L_0$.

On the other hand, we claim that on the event $\calE$,
$$
|\Lambda_L^R|^{-1} \sum_{i\in\Lambda_L^R} \mathbf{1}_{\calB_R^i} \geq 1-2^{-2d-2}.
$$
Indeed, were it not the case, then $|\calA| > 2^{-2d-2}|\Lambda_L^R| = (2/\log\alpha) f(0,\alpha\eta) |\Lambda_L|$, and thus $\calE$ would not occur. Collecting all these estimates, we have proved that
$$
\mu^{+,0}_{2L;0,\eta}(\calB_R^0)\geq 1-2^{-2d-2}.
$$
The conclusion now follows easily. Indeed, FKG implies that pinning all the sites outside $B_R(0)$ only reduces the expectation. Therefore
$$
\bk{\varphi_0}_{L;0,\eta}^{+,0} \geq (1-2^{-2d-2}) \bk{\varphi_0 \given \calB_R^0}_{L;0,\eta}^{+,0} \geq (1-2^{-2d-2}) \bk{\varphi_0}_{2R+1;0,0}^{+,0}.
$$
Now, standard entropic repulsion estimate imply that
$$
\bk{\varphi_0}^{+,0}_{2R+1;0,0} \geq
\begin{cases}
\hat c_2 \log R		&	(d=2),\\
\hat c_d \sqrt{\log R}	&	(d\geq 3),
\end{cases}
$$
for suitable constants $\hat c_d>0$. Of course, $\log R \asymp |\log f(0,\alpha\eta)|$, and the claim is proved.

\bibliographystyle{plain}
\bibliography{wetting-pw-PTRF}

\begin{thebibliography}{10}

\bibitem{Bi2006}
M.~Biskup.
\newblock Reflection positivity and phase transitions in lattice spin models.
\newblock Notes of the lectures given at Prague Summer School on Mathematical
  Statistical Mechanics, September 2006. Arxiv:math-ph/0610025.

\bibitem{BoBr2001}
E.~Bolthausen and D.~Brydges.
\newblock Localization and decay of correlations for a pinned lattice free
  field in dimension two.
\newblock In {\em State of the art in probability and statistics (Leiden,
  1999)}, volume~36 of {\em IMS Lecture Notes Monogr. Ser.}, pages 134--149.
  Inst. Math. Statist., Beachwood, OH, 2001.

\bibitem{BoDeGi2001}
E.~Bolthausen, J.-D. Deuschel, and G.~Giacomin.
\newblock Entropic repulsion and the maximum of the two-dimensional harmonic
  crystal.
\newblock {\em Ann. Probab.}, 29(4):1670--1692, 2001.

\bibitem{BoDeZe1995}
E.~Bolthausen, J.-D. Deuschel, and O.~Zeitouni.
\newblock Entropic repulsion of the lattice free field.
\newblock {\em Comm. Math. Phys.}, 170(2):417--443, 1995.

\bibitem{BoDeZe2000}
E.~Bolthausen, J.-D. Deuschel, and O.~Zeitouni.
\newblock Absence of a wetting transition for a pinned harmonic crystal in
  dimensions three and larger.
\newblock {\em J. Math. Phys.}, 41(3):1211--1223, 2000.
\newblock Probabilistic techniques in equilibrium and nonequilibrium
  statistical physics.

\bibitem{BoVe2001}
E.~Bolthausen and Y.~Velenik.
\newblock Critical behavior of the massless free field at the depinning
  transition.
\newblock {\em Comm. Math. Phys.}, 223(1):161--203, 2001.

\bibitem{CaVe2000}
P.~Caputo and Y.~Velenik.
\newblock A note on wetting transition for gradient fields.
\newblock {\em Stochastic Process. Appl.}, 87(1):107--113, 2000.

\bibitem{De1996}
J.-D. Deuschel.
\newblock Entropic repulsion of the lattice free field. {II}. {T}he
  {$0$}-boundary case.
\newblock {\em Comm. Math. Phys.}, 181(3):647--665, 1996.

\bibitem{DeGiZa2005}
J.-D. Deuschel, G.~Giacomin, and L.~Zambotti.
\newblock Scaling limits of equilibrium wetting models in {$(1+1)$}-dimension.
\newblock {\em Probab. Theory Related Fields}, 132(4):471--500, 2005.

\bibitem{DeVe2000}
J.-D. Deuschel and Y.~Velenik.
\newblock Non-{G}aussian surface pinned by a weak potential.
\newblock {\em Probab. Theory Related Fields}, 116(3):359--377, 2000.

\bibitem{DuMaRiRo1992}
F.~Dunlop, J.~Magnen, V.~Rivasseau, and P.~Roche.
\newblock Pinning of an interface by a weak potential.
\newblock {\em J. Statist. Phys.}, 66(1-2):71--98, 1992.

\bibitem{Fu2005}
T.~Funaki.
\newblock Stochastic interface models.
\newblock In {\em Lectures on probability theory and statistics}, volume 1869
  of {\em Lecture Notes in Math.}, pages 103--274. Springer, Berlin, 2005.

\bibitem{Gi2001}
G.~Giacomin.
\newblock Aspects of statistical mechanics of random surfaces.
\newblock Notes of the lectures given at IHP, fall 2001. Available at\\
  \texttt{http://www.proba.jussieu.fr/pageperso/giacomin/pub/IHP.ps}.

\bibitem{HrVe2004}
O.~Hryniv and Y.~Velenik.
\newblock Universality of critical behaviour in a class of recurrent random
  walks.
\newblock {\em Probab. Theory Related Fields}, 130(2):222--258, 2004.

\bibitem{IoVe2000}
D.~Ioffe and Y.~Velenik.
\newblock A note on the decay of correlations under $\delta$-pinning.
\newblock {\em Probab. Theory Related Fields}, 116(3):379--389, 2000.

\bibitem{La1991}
G.~F. Lawler.
\newblock {\em Intersections of random walks}.
\newblock Probability and its Applications. Birkh\"auser Boston Inc., Boston,
  MA, 1991.

\bibitem{Ve2004}
Y.~Velenik.
\newblock Entropic repulsion of an interface in an external field.
\newblock {\em Probab. Theory Related Fields}, 129(1):83--112, 2004.

\bibitem{Ve2006}
Y.~Velenik.
\newblock Localization and delocalization of random interfaces.
\newblock {\em Probab. Surv.}, 3:112--169 (electronic), 2006.

\end{thebibliography}

\end{document}